\documentclass[a4paper,10pt, reqno]{amsart}
\usepackage{amssymb,latexsym,amsmath,epsfig,amsthm} %% Add other packages as necessary 

\theoremstyle{plain}
\newtheorem{theorem}{Theorem}
\newtheorem*{theorem*}{Theorem}

\newtheorem*{corollary*}{Corollary}
\newtheorem{lemma}{Lemma}
\newtheorem*{lemma*}{Lemma}

\newtheorem*{proposition*}{Proposition}

\newtheorem*{conjecture*}{Conjecture}
\theoremstyle{definition}

\newtheorem*{definition*}{Definition}
\theoremstyle{remark}
\newtheorem{remark}{Remark}
\newtheorem*{remark*}{Remark}
\theoremstyle{example}
\newtheorem{example}{Example}
\newtheorem*{example*}{Example}

\begin{document}
\title[Non-differentiable functions]{Non-differentiable functions defined in terms of classical representations of real numbers}
\author{Symon Serbenyuk}
\address{Institute of Mathematics \\
 National Academy of Sciences of Ukraine \\
  3~Tereschenkivska St. \\
  Kyiv \\
  01004 \\
  Ukraine}
\email{simon6@ukr.net}


% Key words

\begin{abstract}

The present article is devoted to functions from a certain subclass of non-differentiable functions. The arguments and values of considered functions represented by the s-adic representation or the nega-s-adic representation of real numbers. The technique of modeling such functions is the simplest as compared with well-known techniques of modeling non-differentiable functions. In other words, values of these functions are obtained from the s-adic or nega-s-adic representation of the argument by a certain change of digits or combinations of digits.

Integral, fractal  and other properties of considered functions are described.

\end{abstract}
\maketitle

%% \thispagestyle{empty}

%%%%%%%%%%%%%%%%%%%%%%%%%%%%%%%%%%%%%%%%%%%%%%%%%%%%%%%%%%%%%%%%%%%%%%%%

\section{Introduction}

A nowhere differentiable function is a function whose derivative equals infinity or does not exist at each point from the domain of definition.

The idea of  existence of continuous non-differentiable functions appeared in the nineteenth century. In 1854,  Dirichlet speaking at  lectures at Berlin University said on the existence of a continuous function without the derivative.  
In 1830, the first example of a continuous  non-differentiable  function was modeled by Bolzano in  ``Doctrine on Function" but the last paper was published one   hundred years later \cite{{Br1949}, {Kel1955}}.  In 1861, Rieman gives the following example of a non-differentiable  function without proof:
\begin{equation}
\label{eq: the Rieman function}
f(x)=\sum^{\infty} _{n=1}{\frac{sin(n^2x)}{n^2}}.
\end{equation}
Hardy \cite{Hardy1916}, Gerver \cite{Gerver1971}, and Du Bois-Reymond investigated the last-mentioned function. This function has a finite derivative that equals $\frac 1 2$ at points of the form $\xi\pi$, where $\xi$ is a rational number with an odd numerator
and an odd denominator. Function \eqref{eq: the Rieman function} does not have other points of differentiability.

In 1875, Du Bois-Reymond published the following example of function \cite{Reymond1875}: 
$$
f(x)=\sum^{\infty} _{n=1}{a^ncos(b^n\pi x)},
$$
where $0<a<1$ and $b>1$ is an odd integer number such that $ab>1+\frac 3 2 \pi$. The last-mentioned function was modeled by Weierstrass in 1871. This function has the derivative that equals $(+\infty)$  or $(-\infty)$ on an uncountable everywhere dense set.
In the paper \cite{Darboux1875}, own example of non-differentiable function was modeled nearly simultaneously and independently by Darboux 
$$
f(x)=\sum^{\infty} _{n=1}{\frac{sin((n+1)!x)}{n!}}.
$$

In the sequel, other examples of such functions were constructed and classes of non-differentiable functions were founded. 
 The major contribution to these studies was made by the following scientists: Dini {\cite[p. 148--158]{Dini1878}}, Darboux \cite{Darboux1879}, Orlicz \cite{Orlicz1947}, Hankel {\cite[p. 61--65]{Hankel1870}}.

In 1929, the problem on the massiveness of the set of non-differentiable functions in the space of continuous functions was formulated by Steinhaus. In 1931, this problem was solved independently and by different ways  by Banach  \cite{Banach1931} and Mazurkiewicz \cite{Mazurkiewicz1931}. So the following statement is true.  
\begin{theorem}[{Banach-Mazurkiewicz}]
The set of non-differentiable functions in the space $C[0,1]$ of functions, that are continuous on $[0,1]$, with the uniform metric is a set of the second category.
\end{theorem}

There exist also functions  that do not have a finite or infinite one-sided derivative at any point. In 1922, an example of such  function was modeled by Besicovitch  in \cite{Besicovitch1924}. The set of continuous on $[0,1]$ functions whose right-sided derivative equals a finite number  or equals $+\infty$ on an uncountable set is a set of the second Baire category in the space of all continuous functions. Hence the set of functions,  that do not have a finite or infinite one-sided derivative at any point,  is a set of the first category in the space of  continuous on a segment  functions. The last statement was proved by Saks in 1932 (see \cite{Saks1932}).

Now researchers  are trying to find more and more simple examples of non-differentiable functions. Interest in such functions is explained by them connection with fractals, modeling of  real objects, processes, and phenomena (in physics, economics, technology, etc.).

The present article is devoted to simplest examples of non-differentiable functions defined in terms of the s-adic or nega-s-adic representations. 

In addition, let us consider some examples of nowhere differentiable functions defined by another ways.

%%%%%%%%%%%%%%%%%%%%%%%%%%%%%%%%%%%%%%%%%%%%%%%%%

\section{Certain examples of non-differentiable functions}

%%%%%%%%%%%%%%%%%%%%%%%%%%%%%%%%%%%%%%%%%%%%%%%%%

\begin{example}
{\rm Consider the following functions
$$
f\left(\Delta^{3} _{\alpha_1\alpha_2...\alpha_n...}\right)=\Delta^{2} _{\varphi_1(x)\varphi_2(x)...\varphi_n(x)...}
$$
and
$$
g\left(\Delta^{s} _{\alpha_1\alpha_2...\alpha_n...}\right)=\Delta^{2} _{\varphi_1(x)\varphi_2(x)...\varphi_n(x)...},
$$
where $s>2$ is a fixed positive integer number, 
$$
\Delta^{s} _{\alpha_1\alpha_2...\alpha_n...}=\sum^{\infty} _{n=1}{\frac{\alpha_n}{s^n}}, ~\alpha_n\in\{0,1,\dots, s-1\}, 
$$
$$                                                                                          
\varphi_1(x)=\begin{cases}
0 &\text{if $\alpha_1(x)=0$}\\
1&\text{if $\alpha_1(x)\ne0$,}
\end{cases}
$$
and
$$
\varphi_{j}(x)=\begin{cases}
\varphi_{j-1}(x)&\text{for $\alpha_{j}(x)=\alpha_{j-1}(x)$}\\
1-\varphi_{j-1}(x)&\text{for $\alpha_{j}(x)\ne\alpha_{j-1}(x)$.}
\end{cases}
$$
In 1952, the function  $g$ was introduced by Bush in \cite{Bush1952}, and the function $f$ was modeled by Wunderlich in  \cite{Wunderlich1952}.  The functions $f$ and $g$ are non-differentiable.}
\end{example}

In \cite{Salem1943}, Salem modeled the function 
$$
s(x)=s\left(\Delta^2 _{\alpha_1\alpha_2...\alpha_n...}\right)=\beta_{\alpha_1}+ \sum^{\infty} _{n=2} {\left(\beta_{\alpha_n}\prod^{n-1} _{i=1}{q_i}\right)}=y=\Delta^{Q_2} _{\alpha_1\alpha_2...\alpha_n...},
$$
where $q_0>0$, $q_1>0$, and $q_0+q_1=1$. This function is a singular function. However,  
generalizations of the Salem function can be non-differentiable functions or do not have the derivative on a certain set. 

In October 2014,  generalizations of the Salem function such that them arguments represented in terms of positive~\cite{C1869} or alternating \cite{S. Serbenyuk alternating Cantor series 2013} Cantor series   or the nega-$\tilde Q$-representation \cite{{S. Serbenyuk abstract 10},{S16}}  were considered by Serbenyuk in the reports ``Determination of a class of functions represented by  Cantor series by systems of functional equations" and ``Polybasic positive and alternating $\tilde Q$-representations and them applications to determination of functions by systems of  functional equations"  at the fractal analysis seminar of the Institute of Mathematics of the National Academy of Sciences of Ukraine and the National  Pedagogical Dragomanov University (the list of talks available at \\ http://www.imath.kiev.ua/events/index.php?seminarId=21\&archiv=1, the first presentation (in Ukrainian) available at https://www.researchgate.net/publication/314426236). These results were presented by the author in the papers (conference abstracts and articles) \cite{{S. Serbenyuk abstract 9}, {Symon2015}, {Symon2017}, {S. Serbenyuk function nega-tilde Q-representation}} and in the presentation that  available at \\ https://www.researchgate.net/publication/303736670.

Consider these generalizations of the Salem function.
 \begin{example}[\cite{Symon2015}]
{\rm Let $(d_n)$ is a fixed sequence of positive integers, $d_n>1$, and $(A_n)$ is a sequence of the sets  $A_n = \{0,1,\dots,d_n-1\}$.

Let $x\in [0,1]$ be an arbitrary number represented by a positive Cantor series
\begin{equation*}
%\label{series1}
x=\Delta^D _{\varepsilon_1\varepsilon_2...\varepsilon_n...}=\sum^{\infty} _{n=1} {\frac{\varepsilon_n}{d_1d_2\dots d_n}}, ~\mbox{where}~\varepsilon_n \in A_n.
\end{equation*}

Let $P=||p_{i,n}||$  be a fixed matrix such that  $p_{i,n} \ge 0$ ($ n=1,2,\dots ,$ and $i=~\overline{0,d_n-1}$), $\sum^{d_n-1} _{i=0} {p_{i,n}}=1$ for an arbitrary $n \in \mathbb N$, and $\prod^{\infty} _{n=1}{p_{i_n,n}}=0$ for any sequence   $(i_n)$.

Suppose that elements of the matrix $P=||p_{i,n,}||$ can be negative numbers as well but   
$$
\beta_{0,n}=0, \beta_{i,n}>0 ~\mbox{for}~ i\ne 0, ~\mbox{and} ~ \max_i {|p_{i,n}|} <1.
$$
Here 
$$
\beta_{\varepsilon_{k},k}=\begin{cases}
0&\text{if $\varepsilon_{k}=0$}\\
\sum^{\varepsilon_{k}-1} _{i=0} {p_{i,k}}&\text{if $\varepsilon_{k}\ne 0$.}
\end{cases}
$$
Then the following statement is true.
\begin{theorem}
Given the matrix $P$  such that for all $n \in \mathbb N$ the following are true:  $p_{\varepsilon_n,n}\cdot p_{\varepsilon_n-1,n}<0$ moreover $d_n \cdot~p_{d_n-1,n}\ge 1$ or  $d_n \cdot p_{d_n-1,n}\le 1$; and the  conditions 
$$
\lim_{n \to \infty} {\prod^{n} _{k=1} {d_k p_{0,k}}}\ne  0, \lim_{n \to \infty} {\prod^{n} _{k=1} {d_k p_{d_k-1,k}}}\ne 0
$$
hold simultaneously.
Then the function  
$$
F(x)=\beta_{\varepsilon_1(x),1}+\sum^{\infty} _{k=2} {\left(\beta_{\varepsilon_k(x),k}\prod^{k-1} _{n=1} {p_{\varepsilon_n(x),n}}\right)}
$$
is non-differentiable on  $[0,1]$.
\end{theorem}
 }
\end{example}

\begin{example}[\cite{Symon2017}]
{\rm
Let
$P=||p_{i,n}||$ be a given matrix such that  $n=1,2, \dots$ and $i=\overline{0,d_n-1}$. For this matrix the following system of properties  holds: 
$$ 
\left\{
\begin{aligned}
\label{eq: tilde Q 1}
1^{\circ}.~~~~~~~~~~~~~~~~~~~~~~~~~~~~~~~~~~~~~~~~~~~~~~~\forall n \in \mathbb N:  p_{i,n}\in (-1,1)\\
2^{\circ}.  ~~~~~~~~~~~~~~~~~~~~~~~~~~~~~~~~~~~~~~~~~~~~~~~~\forall n \in \mathbb N: \sum^{d_n-1}_{i=0} {p_{i,n}}=1\\
3^{\circ}. ~~~~~~~~~~~~~~~~~~~~~~~~~~~~~~~~~~~~~ \forall (i_n), i_n \in  A_{d_n}: \prod^{\infty} _{n=1} {|p_{i_n,n}|}=0\\
4^{\circ}.~~~~~~~~~~~~~~\forall  i_n \in A_{d_n}\setminus\{0\}: 1>\beta_{i_n,n}=\sum^{i_n-1} _{i=0} {p_{i,n}}>\beta_{0,n}=0.\\
\end{aligned}
\right.
$$ 

 Let us consider the following function 
$$ 
\tilde{F}(x)=\beta_{\varepsilon_1(x),1}+\sum^{\infty} _{n=2} {\left(\tilde{\beta}_{\varepsilon_n(x),n}\prod^{n-1} _{j=1} {\tilde{p}_{\varepsilon_j(x),j}}\right)},
$$
where
$$
\tilde{\beta}_{\varepsilon_n(x),n}=\begin{cases}
\beta_{\varepsilon_n(x),n}&\text{if $n$ is   odd }\\
\beta_{d_n-1-\varepsilon_n(x),n}&\text{if $n$ is  even,}
\end{cases}
$$
$$
\tilde{p}_{\varepsilon_n(x),n}=\begin{cases}
p_{\varepsilon_n(x),n}&\text{if $n$  is odd }\\
p_{d_n-1-\varepsilon_n(x),n}&\text{if $n$  is   even,}
\end{cases}
$$
$$
\beta_{\varepsilon_{n}(x),n}=\begin{cases}
0&\text{if $\varepsilon_{n}=0$}\\
\sum^{\varepsilon_{n}-1} _{i=0} {p_{i,n}}&\text{if $\varepsilon_{n}\ne 0$.}
\end{cases}
$$
Here $x$ represented by an alternating Cantor series, i.e., 
$$
x=\Delta^{-(d_n)} _{\varepsilon_1\varepsilon_2...\varepsilon_n...}=\sum^{\infty} _{n=1} {\frac{1+\varepsilon_n}{d_1d_2\dots d_n}(-1)^{n+1}},
$$
where $(d_n)$ is a fixed sequence of positive integers, $d_n>1$, and $(A_{d_n})$ is a sequence of the sets  $A_{d_n} = \{0,1,\dots,d_n-1\}$, and $\varepsilon_n\in A_{d_n}$.
\begin{theorem}
Let  $p_{\varepsilon_n,n}\cdot p_{\varepsilon_n-1,n}<0$  for all $n \in \mathbb N$, $\varepsilon_n \in A_{d_n} \setminus \{0\}$ and conditions 
$$
\lim_{n \to \infty} {\prod^{n} _{k=1} {d_k p_{0,k}}}\ne  0, \lim_{n \to \infty} {\prod^{n} _{k=1} {d_k p_{d_k-1,k}}}\ne 0
$$
hold simultaneously.  Then the function $\tilde{F}$ is  non-differentiable on $[0,1]$. 
\end{theorem}
 }
\end{example}

\begin{example}[\cite{S. Serbenyuk function nega-tilde Q-representation}]
{\rm Let $\tilde Q=||q_{i,n}||$ be a fixed matrix, where  $i=\overline{0,m_n}$, $m_n \in N^{0} _{\infty}= \mathbb N \cup \{0,\infty\}$, $n=1,2,\dots$, and the following system of poperties is true for elements $q_{i,n}$ of the last-mentioned  matrix
$$
\left\{
\begin{aligned}
\label{eq: tilde Q 1}
1^{\circ}.  ~~~~~~~~~~~~~~~~~~~~~~~~~~~~~~~~q_{i,n}>0\\
2^{\circ}.  ~~~~~~~~~~~~~~~\forall n \in \mathbb N: \sum^{m_n}_{i=0} {q_{i,n}}=1\\
3^{\circ}.  \forall (i_n), i_n \in \mathbb N \cup \{0\}: \prod^{\infty} _{n=1} {q_{i_n,n}}=0.\\
\end{aligned}
\right.
$$

The following expansion  of $x\in[0,1)$
\begin{equation}
\label{def: nega-tilde Q 3}
x= \sum^{i_1-1} _{i=0}{q_{i,1}}+\sum^{\infty} _{n=2}{\left[(-1)^{n-1}\tilde \delta_{i_n,n}\prod^{n-1} _{j=1}{\tilde q_{i_j,j}}\right]}+\sum^{\infty} _{n=1}{\left(\prod^{2n-1} _{j=1}{\tilde q_{i_j,j}}\right)}
\end{equation}
is called \emph{the nega-$\tilde Q$-expansion of $x$}. By $x=\Delta^{-\tilde Q} _{i_1i_1...i_n...}$ denote the nega-$\tilde Q$-expansion of $x$.  The last-mentioned notation is called \emph{the nega-$\tilde Q$-representation of  $x$}. Here
$$
\tilde \delta_{i_{n},n}=\begin{cases}
1&\text{if $n$ is  even   and $i_{n}=m_{n}$}\\
\sum^{m_{n}} _{i=m_{n}-i_{n}} {q_{i,n}}&\text{if $n$ is  even and $i_{n}\ne m_{n}$}\\
0&\text{if $n$ is  odd   and $i_{n}=0$}\\
\sum^{i_n-1} _{i=0}{q_{i,n}}&\text{if $n$ is odd  and $i_{n}\ne0$,}\\
\end{cases}
$$
and the first sum in the expression  \eqref{def: nega-tilde Q 3} is equal to $0$ if  $i_1=0$.

Suppose that $m_n<\infty$ for all positive integers $n$.

Numbers from a some countable subset of $[0,1]$ have two different nega-$\tilde Q$-representations, i.e.,
$$
\Delta^{-\tilde Q} _{i_1i_2...i_{n-1}i_nm_{n+1}0m_{n+3}0m_{n+5}...}=\Delta^{-\tilde Q} _{i_1i_2...i_{n-1}[i_n-1]0m_{n+2}0m_{n+4}...}, ~i_n\ne 0.
$$
These numbers are called  \emph{nega-$\tilde Q$-rationals} and the rest of the numbers from $[0,1]$ are called  \emph{nega-$\tilde Q$-irrationals}.

Let we have matrixes of the same dimension $\tilde Q=||q_{i,n}||$ (properties of the last-mentioned matrix were considered earlier) and  $P=||p_{i,n}||$, where  $i=\overline{0,m_n}$, $m_n\in \mathbb N\cup \{0\}$, $n=1,2,\dots$, and for elements  $p_{i,n}$ of $P$ the following system of conditions is true: 

$$
\left\{
\begin{aligned}
\label{eq: tilde Q 1}
1^{\circ}.~~~~~~~~~~~~~~~~~~~~~~~~~~~~~~~~~~~~~~~~~~~~~~ ~~~~~ p_{i,n}\in (-1,1)\\
2^{\circ}.  ~~~~~~~~~~~~~~~~~~~~~~~~~~~~~~~~~~~~~~~~~~\forall n \in \mathbb N: \sum^{m_n}_{i=0} {p_{i,n}}=1\\
3^{\circ}.  ~~~~~~~~~~~~~~~~~~~~~~~~\forall (i_n), i_n \in \mathbb N \cup \{0\}: \prod^{\infty} _{n=1} {|p_{i_n,n}|}=0\\
4^{\circ}.~~~~~~~~~~~~~~~\forall  i_n \in \mathbb N: 0=\beta_{0,n}<\beta_{i_n,n}=\sum^{i_n-1} _{i=0} {p_{i,n}}<1.\\
\end{aligned}
\right.
$$ 

\begin{theorem}
If the following properties of the matrix  $P$ hold:
\begin{itemize}
\item for all $n \in \mathbb N$, $i_n \in N^1 _{m_n}=\{1,2,\dots,m_n\}$
$$
p_{i_n,n}\cdot p_{i_n-1,n}<0;
$$
\item the conditions 
$$
\lim_{n \to \infty} {\prod^{n} _{k=1} {\frac{p_{0,k}}{q_{0,k}}}}\ne  0, \lim_{n \to \infty} {\prod^{n} _{k=1} {\frac{p_{m_k,k}}{q_{m_k,k}} }}\ne 0
$$
\end{itemize}
hold simultaneously, then the function
$$
F(x)=\beta_{i_1(x),1}+\sum^{\infty} _{k=2}{\left[\tilde \beta_{i_k(x),k}\prod^{k-1} _{j=1}{\tilde p_{i_j(x),j}}\right]}.
$$
 does not have the finite or infinite derivative at any nega-$\tilde Q$-rational point from the segment $[0,1]$.
\end{theorem}
Here $$
\tilde p_{i_n,n}=\begin{cases}
p_{i_n,n}&\text{if $n$ is odd}\\
p_{m_n-i_n,n}&\text{if $n$ is  even, }
\end{cases}
$$
$$
\beta_{i_n,n}=\begin{cases}
\sum^{i_n-1} _{i=0} {p_{i,n}}>0&\text{if $i_n \ne 0$}\\
0&\text{if $i_n=0$,}
\end{cases}~~~
\tilde \beta_{i_n,n}=\begin{cases}
\beta_{i_n,n}&\text{if $n$ is  odd}\\
\beta_{m_n-i_n,n}&\text{if $n$ is even. }
\end{cases}
$$
 }
\end{example}

The last examples of non-differentiable functions are difficult. However, there exist elementary examples of such functions.

%%%%%%%%%%%%%%%%%%%%%%%%%%%%%%%%%%%%%%%%%%%%%%%%%%%%

\section{The simplest example of non-differentiable function, and its analogues}

%%%%%%%%%%%%%%%%%%%%%%%%%%%%%%%%%%%%%%%%%%%%%%%%%%%%

 In 2012, the main results of this section were represented by the author of the present article at  the International Scientific Conference ``Asymptotic Methods in the Theory of Differential Equations" dedicated to 80th anniversary of M.~I.~Shkil \cite{S. Serbenyuk abstract 6} and published in the paper \cite{Symon12(2)} (the working paper available at https://www.researchgate.net/publication/314409844). The version of the last-mentioned published paper into English available at https://arxiv.org/pdf/1703.02820.pdf.

We shall not consider numbers whose  ternary representation has the period $(2)$ (without the number $1$).  Let us consider a certain function  $f$ defined on $[0,1]$ by the following way: 
\begin{equation*}
%\label{ff1}
x=\Delta^{3} _{\alpha_{1}\alpha_{2}...\alpha_{n}...}\stackrel{f}{\rightarrow} \Delta^{3} _{\varphi(\alpha_{1})\varphi(\alpha_{2})...\varphi(\alpha_{n})...}=f(x)=y,
\end{equation*}
where $\varphi(i)=\frac{-3i^{2}+7i}{2}$, $ i \in N^{0} _{2}=\{0,1,2\}$, and $\Delta^{3} _{\alpha_{1}\alpha_{2}...\alpha_{n}...}$ is the ternary  representation of  $x \in [0,1]$. That is values of this function are obtained from the ternary representation
of the argument by the following change of digits: 0 by 0, 1 by 2, and 2 by 1. This function preserves the ternary digit $0$.

In this section,  differential, integral, fractal,  and other properties of the function $f$ are  described; equivalent representations of this function by additionally defined auxiliary functions are considered.

We begin with definitions of some  auxiliary functions.

 Let $i, j, k$ be pairwise distinct digits of the ternary numeral system.
First let us introduce a function $\varphi_{ij} (\alpha)$ defined on the alphabet of the ternary numeral system by the following:
\begin{center}
\begin{tabular}{|c|c|c|c|}
\hline
 &$ i$ &$ j $ & $k$\\
\hline
$\varphi_{ij} (\alpha)$ &$0$ & $0$ & $1$\\
\hline
\end{tabular}
\end{center}
That is  $f_{ij}$  is a function given on $[0;1]$ by the following way
$$
x=\Delta^{3} _{\alpha_{1}\alpha_{2}...\alpha_{n}...}\stackrel{f_{ij}}{\rightarrow} \Delta^{3} _{\varphi_{ij} (\alpha_{1})\varphi_{ij}  (\alpha_{2})...\varphi_{ij} (\alpha_{n})...}=f_{ij} (x)=y.
$$

\begin{remark} From the definition of $f_{ij}$ it follows that 
 $f_{01} =f_{10}$, $f_{02} =f_{20}$, and $f_{12} =f_{21} $. Since it is true, we shall use only the following notations:  $f_{01}, f_{02}, f_{12}$.
\end{remark}

\begin{lemma}
\label{lmf1}
The function $f$ can be represented by the following:
\begin{enumerate}
\item
\begin{equation*}
%\label{ff2} 
f(x)=2x-3f_{01} (x), ~~\mbox{where}~~\Delta^{3} _{\alpha_{1}\alpha_{2}...\alpha_{n}...}\stackrel{f_{01}}{\rightarrow} \Delta^{3} _{\varphi_{01} (\alpha_{1})\varphi_{01} (\alpha_{2})...\varphi_{01} (\alpha_{n})...},
\end{equation*}
$\varphi_{01} (i)=\frac{i^2 - i}{2}$, $i \in N^0 _2$;
\item
\begin{equation*}
%\label{ff3}
f(x)=\frac{3}{2}-x-3f_{12} (x), ~~\mbox{where}~~\Delta^{3} _{\alpha_{1}\alpha_{2}...\alpha_{n}...}\stackrel{f_{12}}{\rightarrow} \Delta^{3} _{\varphi_{12} (\alpha_{1})\varphi_{12} (\alpha_{2})...\varphi_{12} (\alpha_{n})...},
\end{equation*}
$\varphi_{12} (i)=\frac{i^2 - 3i+2}{2}$, $i \in N^0 _2$.

\item
\begin{equation*}
%\label{ff4}
f(x)=\frac{x}{2}+\frac{3}{2}f_{02} (x), ~~\mbox{where}~~\Delta^{3} _{\alpha_{1}\alpha_{2}...\alpha_{n}...}\stackrel{f_{02}}{\rightarrow} \Delta^{3} _{\varphi_{02} (\alpha_{1})\varphi_{02} (\alpha_{2})...\varphi_{02} (\alpha_{n})...},
\end{equation*}
$\varphi_{02} (i)=-i^2 +2i$, $i \in N^0 _2$.
\end{enumerate}
\end{lemma}

\begin{lemma}
\label{lmf2}
The functions $f, f_{01},f_{02}, f_{12}$ have the following properties:
\begin{enumerate}
\item
$$
[0,1]\stackrel{f}{\rightarrow} \left([0,1] \setminus \{\Delta^3 _{\alpha_1\alpha_2...\alpha_n111...}\}\right) \cup\left\{\frac{1}{2}\right\};
$$
\item the point $x_0=0$ is the unique invariant point of the function  $f$;

\item the function $f$ is not bijective on a certain countable subset of $[0,1]$.
\item the following relationships hold for all $x \in [0,1]$:
\begin{equation*}
%\label{ff5}
f(x)-f(1-x)=f_{01} (x)-f_{12} (x),
\end{equation*}
\begin{equation*}
%\label{ff6}
f(x)+f(1-x)=\frac{1}{2}+3f_{02} (x),
\end{equation*}
\begin{equation*}
%\label{ff7}
f_{01} (x)+f_{02} (x)+f_{12} (x)=\frac{1}{2},
\end{equation*}
\begin{equation*}
%\label{ff8}
2f_{01} (x)+f_{02} (x)=x,
\end{equation*}
\begin{equation*}
 %\label{ff9}
f_{01} (x)- f_{12} (x)=x- \frac{1}{2};
\end{equation*}

\item the function $f$ is not monotonic on the domain of definition; in particular, the function $f$ is a decreasing function on the set
$$
\{x: x_1<x_2 \Rightarrow (x_1=\Delta^3 _{c_1...c_{n_0}1\alpha_{n_0+2}\alpha_{n_0+3}...} \wedge  x_2=\Delta^3 _{c_1...c_{n_0}2\beta_{n_0+2}\beta_{n_0+3}...})\}, 
$$
 where $n_0 \in \mathbb Z_0=\mathbb N \cup \{0\}$, $c_1, c_2,\dots, c_{n_0}$ is an ordered set of the ternary digits, $\alpha_{n_0+i} \in N^0 _2$,  $\beta_{n_0+i} \in N^0 _2$, $i \in \mathbb N$;
and the function $f$ is an increasing function on the set
$$
\{x: x_1<x_2 \Rightarrow (x_1=\Delta^3 _{c_1...c_{n_0}0\alpha_{n_0+2}\alpha_{n_0+3}...} \wedge  x_2=\Delta^3 _{c_1...c_{n_0}r\beta_{n_0+2}\beta_{n_0+3}...})\}, 
$$
where $r \in \{1,2\}$.
\end{enumerate}
\end{lemma}

Let us consider fractal properties of all level sets of the functions $f_{01}, f_{02}, f_{12}$. 

The following set  
$$
f^{-1} (y_0)=\{x: g(x)=y_0\},
$$
where $y_0$ is a fixed element of the range of values $E(g)$ of the function $g$, is called \emph{a level set of  $g$}.

\begin{theorem}
\label{rivnif1}
The following statements are true:
\begin{itemize}
\item if there exists at least one digit $2$ in the ternary representation of $y_0$, then $f^{-1} _{ij} (y_0)=\varnothing$;

\item if $y_0=0$ or $y_0$  is a ternary-rational number from the set $C[3, \{0,1\}]$, then  
$$
\alpha_0(f^{-1} _{ij} (y_0))=\log_3 2;
$$

\item if $y_0$ is a ternary-irrational number from the set $C[3, \{0,1\}]$, then
$$
0 \le \alpha_0(f^{-1} _{ij} (y_0))\le \log_3 2,
$$
where  $\alpha_0(f^{-1} _{ij} (y_0))$ is the Hausdorff-Besicovitch  dimension of  $f^{-1} _{ij} (y_0)$.
\end{itemize}
\end{theorem}

Let us describe the main properties of the function $f$.
\begin{theorem}
The function $f$ is continuous at ternary-irrational points, and ternary-rational points are points of discontinuity of the function.
Furthermore, a ternary-rational point $x_0=\Delta^3 _{\alpha_1\alpha_2...\alpha_n000...} $ is a  point of discontinuity $\frac{1}{2\cdot 3^{n-1}}$ whenever $\alpha_n=1$, and is a point of discontinuity $\left(-\frac{1}{2\cdot3^{n-1}}\right)$ whenever $\alpha_n=2$.
\end{theorem}

\begin{theorem}
The function $f$ is non-differentiable. 
\end{theorem}

Let us consider one fractal property of the graph of $f$.
Suppose that 
$$
X=[0;1]\times[0;1]=\left\{(x,y): x=\sum^{\infty} _{m=1} {\frac{\alpha_m}{3^{m}}}, \alpha_{m} \in N^{0} _{2},
y=\sum^{\infty} _{m=1} {\frac{\beta_m}{3^{m}}}, \beta_{m} \in N^{0} _{2}\right\}.
$$
Then the set 
$$
\sqcap_{(\alpha_{1}\beta_{1})(\alpha_{2}\beta_{2})...(\alpha_{m}\beta_{m})}=\Delta^{3} _{\alpha_{1}\alpha_{2}...\alpha_{m}}\times\Delta^{3} _{\beta_{1}\beta_{2}...\beta_{m}}
$$
is a square with a side length of $3^{-m}$. This square is called \emph{a square of rank $m$ with  base $(\alpha_{1}\beta_{1})(\alpha_{2}\beta_{2})\ldots (\alpha_{m}\beta_{m})$}.

If  $E\subset X$, then the number
$$
\alpha^{K}(E)=\inf\{\alpha: \widehat{H}_{\alpha} (E)=0\}=\sup\{\alpha: \widehat{H}_{\alpha} (E)=\infty\},
$$
where
$$
\widehat{H}_{\alpha} (E)=\lim_{\varepsilon \to 0} \left[{\inf_{d\leq \varepsilon} {K(E,d)d^{\alpha}}}\right]
$$
and $K(E,d)$ is the minimum number of squares of diameter $d$ required to cover the set $E$, is called \emph{ the fractal cell entropy dimension of the set E.} It is easy to see that $\alpha^{K}(E)\ge \alpha_0(E)$.

The notion of  the fractal cell entropy dimension  is used for the calculation of the Hausdorff-Besicovitch dimension of the graph of $f$, because, in the case of the function $f$, we obtain that $\alpha^{K}(E)= \alpha_0(E)$ (it follows from the self-similarity of the graph  of $f$). 
\begin{theorem}
The Hausdorff-Besicovitch dimension of the graph of $f$  is equal to 1.
\end{theorem}

Integral properties of $f$ is described in the following theorem.
\begin{theorem}
The Lebesgue integral of the function $f$   is equal to $\frac{1}{2}$.
\end{theorem}

There exist several analogues of the function $f$ such that have the same properties as this function and are defined by analogy. 
Let us consider these functions.

One can define $m=3!=6$ functions  determined on $[0,1]$ in terms of the ternary numeral system by the following way:
$$
\Delta^{3} _{\alpha_{1}\alpha_{2}...\alpha_{n}...}\stackrel{f_m}{\rightarrow} \Delta^{3} _{\varphi_m(\alpha_{1})\varphi_m(\alpha_{2})...\varphi_m(\alpha_{n})...},
$$
where the function $\varphi_m(\alpha_n)$ determined on an alphabet of the ternary numeral system and $f_m (x)$ is defined by the following table for each $m=\overline{1,6}$. 
\begin{center}
\begin{tabular}{|c|c|c|c|}
\hline
 &$ $ 0 &$ 1 $ & $2$\\
\hline
$\varphi_1 (\alpha_n) $ &$0$ & $1$ & $2$\\
\hline
$\varphi_2 (\alpha_n) $ &$0$ & $2$ & $1$\\
\hline
$\varphi_3 (\alpha_n) $ &$1$ & $0$ & $2$\\
\hline
$\varphi_4 (\alpha_n) $ &$1$ & $2$ & $0$\\
\hline
$\varphi_5 (\alpha_n) $ &$2$ & $0$ & $1$\\
\hline
$\varphi_6 (\alpha_n) $ &$2$ & $1$ & $0$\\
\hline
\end{tabular}
\end{center}
That is one can model a class of functions whose values  are obtained from the ternary representation
of the argument by  a certain  change of ternary digits.

It is easy to see that the function $f_1 (x)$ is the function  $y=x$ and the function $f_6 (x)$ is the function $y=1-x$, i.e.,
$$
y=f_1(x)=f_1\left(\Delta^3 _{\alpha_1\alpha_2\ldots \alpha_n\ldots}\right)=\Delta^3 _{\alpha_1\alpha_2\ldots \alpha_n\ldots}=x,
$$
$$
y=f_6(x)=f_6\left(\Delta^3 _{\alpha_1\alpha_2\ldots \alpha_n\ldots}\right)=\Delta^3 _{[2-\alpha_1][2-\alpha_2]\ldots  [2-\alpha_n]\ldots}=1-x.
$$

We shall describe  some application of function of the last form in the next subsection.

\begin{lemma} 
Any function  $f_m$ can be represented by the functions $f_{ij}$ in the following form
$$
f_m=a^{(ij)} _{m}x+b^{(ij)} _{m}+c^{(ij)} _{m}f_{ij} (x), ~\mbox{where}~a^{(ij)} _{m}, b^{(ij)} _{m}, c^{(ij)} _{m} \in \mathbb Q.
$$
\end{lemma}

One can  formulate the following corollary.

\begin{theorem}
The function $f_m $ such that $f_m(x)\ne x$ and $f_m(x)\ne1-x$ is:
\begin{itemize}
\item continuous almost everywhere;
\item non-differentiable on $[0,1]$;
\item a function whose the Hausdorff-Besicovitch dimension of the grapf is equal to  $1$;
\item  a function whose the Lebesgue integral is equal to $\frac{1}{2}$.
\end{itemize}
\end{theorem}

Generalizations of  results described in this section will be considered in the following subsection.

%%%%%%%%%%%%%%%%%%%%%%%%%%%%%%%%%%%%%%%%%%%%%%

\section{Generalizations of the simplest example of non-differentiable function}

%%%%%%%%%%%%%%%%%%%%%%%%%%%%%%%%%%%%%%%%%%%%%%

 In 2013, the investigations of the  last section were generalized by the author in  several conference abstracts \cite{{S. Serbenyuk abstract 7}, {S. Serbenyuk abstract 8}} and  in the paper \cite{S. Serbenyuk functions with complicated local structure 2013} ``One one class of functions with complicated local structure". Consider these results. 

We begin with definitions.

Let $s>1$ be a fixed positive integer number, and let the set $A=\{0,1,\dots,s-~1\}$ be an alphabet of the s-adic or nega-s-adic numeral system. The notation $x=\Delta^{\pm s} _{\alpha_1\alpha_2...\alpha_n...}$ means that $x$ is represented by the s-adic or nega-s-adic representation, i.e.,
$$
x=\sum^{\infty} _{n=1}{\frac{\alpha_n}{s^n}}=\Delta^{ s} _{\alpha_1\alpha_2...\alpha_n...}
$$
or
$$
x=\sum^{\infty} _{n=1}{\frac{(-1)^n\alpha_n}{s^n}}=\Delta^{ -s} _{\alpha_1\alpha_2...\alpha_n...}, \alpha_n\in A.
$$

Let $\Lambda_{s}$ be a class of functions of the type 
\begin{equation*}
  %\label{eq:1}
 f: x=\Delta^{\pm s} _{\alpha_1\alpha_2...\alpha_n...} \to \Delta^{\pm s} _{\beta_1\beta_2...\beta_n...}=f(x)=y,
\end{equation*}
where $\left(\beta_{km+1},\beta_{km+2},\dots,\beta_{(m+1)k}\right)=\theta\left(\alpha_{km+1},\alpha_{km+2},\dots,\alpha_{(m+1)k}\right)$, the number $k$ is a fixed positive integer for a specific function $f$, $m=0,1,2,\dots,$ and
$\theta(\gamma_1,\gamma_2,\dots,\gamma_k)$ is some function of $k$ variables (it is the bijective correspondence) such that the  set  
$$
A^k=\underbrace{A \times A \times \ldots \times A}_{k}.
$$
is its domain of definition and range of values.

Each combination $(\gamma_1,\gamma_2,\dots,\gamma_k)$ of $k$ s-adic or nega-s-adic digits (according to the number representation of the argument of a function $f$) is assigned to the single combination $\theta(\gamma_1,\gamma_2,\dots,\gamma_k)$ of $k$ s-adic or nega-s-adic digits  (according to the number representation of the value of a function $f$). The combination $\theta(\gamma_1,\gamma_2,\dots,\gamma_k)$ is assigned to the unique combination $(\gamma^{'} _1,\gamma^{'} _2,\dots,\gamma^{'} _k)$ that may be not to match with 
$(\gamma_1,\gamma_2,\dots,\gamma_k)$. The $\theta$ is a bijective function on $A^k$.

 It is clear that any function  $f \in \Lambda_{s}$ is  one of the following functions:
$$
f^s _k, ~~~f_+,~~~ f^{-1} _+,~~~ f_+ \circ f^s _k,~~~ f^s _k \circ f^{-1} _+, ~~~ f_+ \circ f^s _k \circ  f^{-1} _+, 
$$
where
\begin{equation*}
%\label{form: f^s _k1}
f^s _k\left(\Delta^s _{\alpha_1\alpha_2...\alpha_n...}\right)=  \Delta^{s}_{\beta_{1}\beta_{2}...\beta_{n}...}, 
\end{equation*}
$\left(\beta_{km+1},\beta_{km+2},\dots,\beta_{(m+1)k}\right)=\theta\left(\alpha_{km+1},\alpha_{km+2},\dots,\alpha_{(m+1)k}\right)$ for $m=0,1,2,\dots,$ and  some fixed positive integer number $k$, i.e., 
$$
\left(\beta_{1},\beta_{2},\dots ,\beta_{k}\right)=\theta\left(\alpha_{1},\alpha_{2},\dots ,\alpha_{k}\right),
$$
$$
\left(\beta_{k+1},\beta_{k+2},\dots ,\beta_{2k}\right)=\theta\left(\alpha_{k+1},\alpha_{k+2},\dots ,\alpha_{2k}\right),
$$
$$
\dots \dots \dots \dots \dots \dots \dots 
$$
$$
\left(\beta_{km+1},\beta_{km+2},\dots ,\beta_{(m+1)k}\right)=\theta\left(\alpha_{km+1},\alpha_{km+2},\dots ,\alpha_{(m+1)k}\right),
$$
$$
\dots \dots \dots \dots \dots \dots \dots 
$$
and
 \begin{equation*}
%\label{form: f _+1}
f_+\left(\Delta^s _{\alpha_1\alpha_2...\alpha_n...}\right)=\Delta^{-s} _{\alpha_1\alpha_2...\alpha_n...},
\end{equation*}
\begin{equation*}
%\label{form: f^{-1} _+1}
f^{-1} _+\left(\Delta^{-s} _{\alpha_1\alpha_2...\alpha_n...}\right)= \Delta^s _{\alpha_1\alpha_2...\alpha_n...}.
\end{equation*}

Let us consider several examples.

The function $f$ considered in the last subsection  is a function  of the $f^3 _1$ type. In fact, 
$$
x=\Delta^3 _{\alpha_1\alpha_2...\alpha_n...} \stackrel{f}{\rightarrow} \Delta^{3}_{\varphi(\alpha_{1})\varphi(\alpha_{2})...\varphi(\alpha_{n})...}=f(x)=y,
$$
where $\varphi\left(\alpha_{n}\right)$ is a function defined in terms of the s-adic numeral system in the following way:

\begin{center}
\begin{tabular}{|c|c|c|c|}
\hline
$\alpha_n$ & $ 0 $& $1$ & $2$\\
\hline
$\varphi (\alpha_n)$ & $0$ & $2$ & $1$\\
\hline
\end{tabular}
\end{center}

Now, we present the example of the function $f^2 _2$. The following function
$$
f^2 _2: \Delta^2 _{\alpha_1\alpha_2...\alpha_n...} \rightarrow \Delta^2 _{\beta_1\beta_2...\beta_n...},
$$ 
where $(\beta_{2m+1},\beta_{2(m+1)})=\theta(\alpha_{2m+1},\alpha_{2(m+1)})$, $m=0,1,2,3,\dots,$ and 
\begin{center}
\begin{tabular}{|c|c|c|c|c|c|}
\hline
$\alpha_{2m+1}\alpha_{2(m+1)}$ & $ 00 $& $01$ & $10$ & $11$\\
\hline
$\beta_{2m+1}\beta_{2(m+1)}$ & $10$ & $11$ & $00$ & $01$\\
\hline
\end{tabular}
\end{center}
 is an example of the $f^2 _2$-type function.

It is obvious that the set of $f^2 _1$ functions consists only of the functions $y=x$ and $y=1-x$ in the binary numeral system. But the set of $f^2 _2$ functions has the order, that is equal to   $4!$, and includes the functions $y=x$ and  $y=1-x$ as well.

\begin{remark}
{\rm The class  $\Lambda_s$ of functions includes the following  linear functions: $y=x$, 
$$
y=f(x)=f\left(\Delta^s _{\alpha_1\alpha_2...\alpha_n...}\right)=\Delta^s _{[s-1-\alpha_1][s-1-\alpha_2]...[s-1-\alpha_n]...}=1-x,
$$
$$
y=f(x)=f\left(\Delta^{-s} _{\alpha_1\alpha_2...\alpha_n...}\right)= \Delta^{-s} _{[s-1-\alpha_1][s-1-\alpha_2]...[s-1-\alpha_n]...}=-\frac{s-1}{s+1}-x.
$$
These functions are called \emph{$\Lambda_s$-linear functions.}}
\end{remark}

\begin{remark}
{\rm The last two functions in the last-mentioned remark are interesting for applications in certain   investigations. For example, in the case of a positive Cantor series, such type function  is of the form:
\begin{equation*}
\begin{split}
f(x) &=f\left(\Delta^D _{\varepsilon_1\varepsilon_2...\varepsilon_n...}\right)=f\left(\sum^{\infty} _{n=1}{\frac{\varepsilon_n}{d_1d_2\dots d_n}}\right)\\
&=\Delta^D _{[d_1-1-\varepsilon_1][d_2-1-\varepsilon_2]...[d_n-1-\varepsilon_n]...}=\sum^{\infty} _{n=1}{\frac{d_n-1-\varepsilon_n}{d_1d_2\dots d_n}}.
\end{split}
\end{equation*}
It is easy to see the last function is a transformation preserving the Hausdorff-Besicovitch dimension.

Consider the following representations by alternating Cantor series:
$$
\Delta^{-D} _{\varepsilon_1\varepsilon_2...\varepsilon_n...}=\sum^{\infty} _{n=1}{\frac{(-1)^n\varepsilon_n}{d_1d_2\dots d_n}},
$$
$$
\Delta^{-(d_n)} _{\varepsilon_1\varepsilon_2...\varepsilon_n...}=\sum^{\infty} _{n=1}{\frac{1+\varepsilon_n}{d_1d_2\dots d_n}(-1)^{n+1}}.
$$

In September 2013, the investigation of relations between  positive   and  alternating   Cantor series (and another investigations of alternating Cantor series) were presented by the author of the present article at the fractal analysis seminar of the Institute of Mathematics of the National Academy of Sciences of Ukraine and the National  Pedagogical Dragomanov University  (see paper \cite{S. Serbenyuk alternating Cantor series 2013}, the presentation and the working paper available at https://www.researchgate.net/publication/303720347 and https://www.researchgate.net/publication/316787375  respectively). 

Consider the following results that follow from relations between  positive   and  alternating   Cantor series.
\begin{lemma}
The following functions are identity transformations:
\end{lemma}
$$
x=\Delta^{D} _{\varepsilon_1\varepsilon_2\ldots\varepsilon_n\ldots}\stackrel{f}{\rightarrow}\Delta^{-(d_n)} _{\varepsilon_1[d_2-1-\varepsilon_2]\ldots\varepsilon_{2n-1}[d_{2n}-1-\varepsilon_{2n}]\ldots}=f(x)=y,
$$
$$
x=\Delta^{-(d_n)} _{\varepsilon_1\varepsilon_2\ldots\varepsilon_n\ldots}\stackrel{g}{\rightarrow}\Delta^{D} _{\varepsilon_1[d_2-1-\varepsilon_2]\ldots\varepsilon_{2n-1}[d_{2n}-1-\varepsilon_{2n}]\ldots}=g(x)=y.
$$

Therefore the  following functions are DP-functions (functions preserving the fractal Hausdorff-Besicovitch dimension): 
$$
x=\Delta^{D} _{\varepsilon_1\varepsilon_2\ldots\varepsilon_n\ldots}\stackrel{f}{\rightarrow}\Delta^{-(d_n)} _{[d_1-1-\varepsilon_1]\varepsilon_2\ldots [d_{2n-1}-1-\varepsilon_{2n-1}]\varepsilon_{2n}\ldots}=f(x)=y,
$$
$$
x=\Delta^{-(d_n)} _{\varepsilon_1\varepsilon_2\ldots\varepsilon_n\ldots}\stackrel{g}{\rightarrow}\Delta^{D} _{[d_1-1-\varepsilon_1]\varepsilon_2\ldots[d_{2n-1}-1-\varepsilon_{2n-1}]\varepsilon_{2n}\ldots}=g(x)=y.
$$
}
\end{remark}

A new method for the construction of metric, probabilistic and dimensional theories for families of representations of real numbers via studies of special mappings (G-isomorphisms of representations), under which symbols of a given representation are mapped into the same symbols of other representation from the same family, and they preserve the Lebesgue measure and the Hausdorff-Besicovitch dimension follows from  investigations considered in this remark  and investigations of functions $f_+, f^{-1} _+$.

Let us describe the main properties of functions $f\in \Lambda_{s}$.
  
\begin{lemma}
For any function  $f$  from $\Lambda_{s}$ except for $\Lambda_s$-linear functions, values of function $f$ for different representations of s-adic rational numbers from ~$[0,1]$ (nega-s-adic rational numbers from $[-\frac{s}{s+1},\frac{1}{s+1}]$ respectively) are different.
\end{lemma}

\begin{remark}
{\rm From the unique representation for each s-adic irrational number from $[0,1]$ it follows that the  function $f^s_k$ is well defined at s-adic irrational  points.

To reach that any function $f \in \Lambda_{s}$ such that $f(x)\ne x$ and $f(x)\ne1-x$ be well-defined on the set of s-adic rational numbers from $[0,1]$, we shall not consider the s-adic representation, which has period $(s-1)$.

Analogously, we shall not consider the nega-s-adic representation, which has period $(0[s-1])$.}
\end{remark}

\begin{lemma}
The set of functions  $f^s _k$ with the defined operation ``composition of functions'' is a finite group that has order equal to $\left(s^k\right)!$.
\end{lemma}

\begin{lemma}
The function $f \in \Lambda_{s}$ such that $f(x)\ne x$, $f(x)\ne -\frac{s-1}{s+1}-x$, and $f(x)\ne 1-x$ has the following properties:
\begin{enumerate}
\item $f$ reflects  $[0,1]$ or $[-\frac{s}{s+1},\frac{1}{s+1}]$ (according to the number representation of the argument of a function $f$) into one of the segments $[0,1]$ or $[-\frac{s}{s+1},\frac{1}{s+1}]$  without enumerable subset of points (according to the number representation of the value of a function $f$) .
\item the function $f$ is not monotonic on the domain of definition;
\item the function $f$ is not a bijective mapping on the domain of definition.
\end{enumerate}
\end{lemma}

\begin{lemma}
The following properties of the  set of invariant points of the function $f^s_k$ are true:
\begin{itemize}
\item the set of invariant points of $f^s _k$ is a continuum set, and its Hausdorff--Besicovitch dimension is equal to  $\frac{1}{k}\log_s j$, when there exists the set  $\{\sigma_1,\sigma_2,...,\sigma_j\}$ $(j \ge 2)$ of k-digit combinations $\sigma_1,...,\sigma_j$ of s-adic digits such that
$$
\theta(a^{(i)} _1,a^{(i)} _2,...,a^{(i)} _k)=(a^{(i)} _1,a^{(i)} _2,...,a^{(i)} _k),~\mbox{where} ~\sigma_i=(a^{(i)} _1 a^{(i)} _2...a^{(i)} _k),~i=\overline{1,j};
$$

\item the  set of invariant points of $f^s_k$ is a finite set, when there exists the unique k-digit combination $\sigma$ of s-adic digits such that
$$
\theta(a_1,a_2,...,a_k)=(a_1,a_2,...,a_k),~\sigma=(a_1a_2...a_k);
$$

\item the  set of invariant points of $f^s_k$ is an empty set, when there not exist any k-digit combination $\sigma$ of s-adic digits such that
$$
\theta(a_1,a_2,...,a_k)=(a_1,a_2,...,a_k),~\sigma=(a_1a_2...a_k).
$$
\end{itemize}
\end{lemma}

In addition, the functions $ f_+$ and $ f^{-1} _+$ have the following properties.
\begin{lemma}
\label{lm4}
 For each $x \in [0,1]$, the function $ f_+$ satisfies the  equation
\begin{equation*}
%\label{eq1}
f(x) +f (1-x)=-\frac{s-1}{s+1};
\end{equation*}
\end{lemma}

\begin{lemma}
For each $y \in [-\frac{s}{s+1},\frac{1}{s+1} ]$, the function   $ f^{-1} _+$ satisfies the equation
\begin{equation*}
%\label{eq2}
f^{-1} (y) +f^{-1}  \left(-\frac{s-1}{s+1}-y\right)=1.
\end{equation*}
\end{lemma}

\begin{lemma}
The set of invariant points of the function $f_{+},$ as well as
$f^{-1} _{+},$ is a self-similar fractal, and its Hausdorff--Besicovitch
dimension is equal to $\frac1 2$.
\end{lemma}

The following theorems are the main theorems about    properties  of functions $f\in \Lambda_s$.

\begin{theorem}
A function   $f \in \Lambda_{s}$ such that $f(x)\ne x$, $f(x)\ne -\frac{s-1}{s+1}-x$, and $f(x)\ne 1-x$  is:
\begin{itemize}
\item continuous at s-adic irrational or nega-s-adic irrational points, and s-adic rational or nega-s-adic rational points are points of discontinuity of this function (according to the number representation of the argument of the function $f$);

\item a non-differentiable function.
\end{itemize}
\end{theorem}

\begin{theorem} Let  $f \in \Lambda_{s}$. Then the following are true: 
\begin{itemize}
\item the Hausdorff--Besicovitch dimension of the graph of any function from the class  $\Lambda_{s}$ is equal to~$1$;
\item
 $$
\int\limits_{D(f)} f(x) \, \mathrm dx=\frac{1}{2},~\mbox{where  $D(f)$ is the domain of definition of $f$.}
$$
\end{itemize}
\end{theorem}

So, in the present article, we considered historical moments of the development of the theory of non-differentiable functions, difficult and simplest examples of such functions. Integral, fractal, and other properties of the simplest example of  nowhere differentiable function and its analogues and generalizations are described. Equivalent representations of the considered simplest example by additionally defined auxiliary functions were reviewed.

\end{document}